\documentclass[10pt,leqno]{amsart}

\newtheorem{lemma}{Lemma}
\newtheorem{theorem}{Theorem}
\newtheorem{corollary}{Corollary}
\newtheorem{proposition}{Proposition}
\newtheorem{definition}{Definition}

\usepackage{graphics}
\usepackage{epsfig}
\usepackage{amssymb,amsfonts,amscd}
\usepackage{enumerate}
\def\C{{\mathbb C}}

\def\CC{{\mathcal C}}
\begin{document}


\title[Wong-Rosay Theorem in almost complex manifolds]
{Wong-Rosay Theorem in almost complex manifolds}

\author{Herv\'e Gaussier and Alexandre Sukhov}

\address{}

\email{} 


\date{\number\year-\number\month-\number\day}

\begin{abstract}
We study the compactness of sequences of diffeomorphisms in almost complex
manifolds in terms of the direct images of the standard integrable structure.
\end{abstract}

\maketitle
\section*{Introduction}
The classical Wong-Rosay theorem states that
every domain in the euclidean space $\C^n$ or more generally in a complex 
manifold of dimension $n$, with an
automorphism orbit accumulating at a $\CC^2$ strictly pseudoconvex
point, is biholomorphically equivalent to the unit ball in $\C^n$ 
(see \cite{gkk02, pi89, Ro, Wo}). 
The aim of this paper is to extend this theorem to strictly
pseudoconvex domains in almost complex
manifolds. Our main result can be considered as a compactness theorem for
sequences of diffeomorphisms. This shows that the convergence of
such sequences can be controlled in terms of the direct images of the standard
complex structure. Our approach is based on the 
scaling method, introduced by S.Pinchuk \cite{pi89} for the case of the 
integrable structure. In order\ to apply it in the almost complex case, we
need substantial modifications. 
In particular we use lower estimates for the Kobayashi
infinitesimal pseudometric on almost complex manifolds \cite{gs03}
and a priori estimates
for $J$-holomorphic curves \cite{sik94}.

\section{Statement of the results}
An almost complex manifold 
$(M,J)$ is a smooth ($\mathcal C^\infty$) real manifold equipped with
an almost complex structure $J$, that is a $\mathcal C^\infty$-field
of complex linear structures on the tangent bundle $TM$ of $ M$.
Given two almost complex manifolds $(M,J)$ and $(M',J')$ and a smooth map $f$ 
from $M'$ to $M$ we say that $f$ is {\sl
$(J',J)$-holomorphic} if its differential $df : TM' \rightarrow TM$
satisfies $df \circ J' = J \circ df$ on $TM$. 
We denote by $\mathcal O_{(J',J)}(M',M)$ the set of
$(J',J)$-holomorphic maps from $M'$ to $M$ and by
$Diff_{(J',J)}(M',M)$ the set of $(J',J)$-holomorphic diffeomorphisms
from $M'$ to $M$. The set $Diff_{(J',J)}(M',M)$ is generically empty. However
given a diffeomorphism $f$ from $M'$ to $M$ and an almost complex
structure $J'$ on $M'$ then $f \in Diff_{(J,J')}(M,M')$ for the almost
complex structure $J=df \circ J' \circ df^{-1}$ naturally associated
with $f$. If $M$ and $M'$ are two domains in $\C^n$, if $J=J'=J_{st}$,
the usual standard structure on $\C^n$, and $f \in
Diff_{(J_{st},J_{st})}(M',M)$ we simply say that $f$ is a biholomorphism
from $M'$ to $M$, or that $M$ is biholomorphic to $M'$.

\vskip 0,1cm
We consider the following situation :

\noindent $\bullet$ $D$ is a bounded domain in $\C^2$,

\noindent $\bullet$ $p$ is a point in a four dimensional almost complex
manifold $(M,J)$, $U$ is a relatively compact neighborhood of
$p$ in $M$, $r$ is a $\mathcal C^2$ strictly 
$J$-plurisubharmonic function in a neighborhood
of $\bar{U}$ satisfying $r(p) =0$ and $dr \neq 0$ on $U$,

\noindent $\bullet$ $(r_\nu)_\nu$ is a sequence of $\mathcal C^2$ functions
in a neighborhood of $\bar{U}$ such that 
$\lim_{\nu \rightarrow \infty}\|r_\nu -r\|_{\mathcal C^2(\bar{U})} =0$,

\noindent $\bullet$ for every $\nu$, $f^\nu$ is a diffeomorphism from
$D$ to $G^\nu \subset M$ with 
$G^\nu \cap U =\{x \in U : r_\nu(x) < 0\}$

\noindent $\bullet$ for every $\nu$, the almost complex structure 
$J_\nu :=df^\nu \circ J_{st} \circ (df^\nu)^{-1}$ extends smoothly to
$\bar{U}$.

\vskip 0,1cm
Then we have :

\begin{theorem}\label{thm2.1}
Assume that $\lim_{\nu \rightarrow \infty}
\|J_\nu - J\|_{\mathcal C^2(\bar{U})} = 0$. 
If there is a point $x^0 \in D$ such that $\lim_{\nu
\rightarrow \infty}f^\nu(x^0) = p$, then $D$ is biholomorphic to the
unit ball $\mathbb B_2$ in $\C^2$.
\end{theorem}


\vskip 0,1cm
Theorem~\ref{thm2.1} implies the following compactness
result : 
\begin{corollary}\label{cor2.2} Let $D$ be a bounded domain in $\C^2$, not
biholomorphic to the unit ball $\mathbb B_2$ and let 
$G=\{x \in M : r(x)<0\}$ be a relatively compact domain in $(M,J)$, 
where $r$ is a $\CC^2$ strictly $J$-plurisubhamonic function on $M$. 
Let $(f^\nu)_\nu$ be a 
sequence of diffeomorphisms from $D$ to $G$, defined in a neighborhood 
of $\bar{D}$. If $\|J_\nu-J\|_{\CC^2(\bar{D})} \rightarrow_{\nu 
\rightarrow \infty}
 0$ then $(f^\nu)_\nu$ is compact.
\end{corollary}
Finally we have an analogue of the Wong-Rosay theorem in
almost complex manifolds :
\begin{theorem}\label{thm2.3}
Let $D=\{x \in M : r(x) < 0\}$ be a relatively compact domain in a four
dimensional almost complex manifold $(M,J)$, where $r$ is a $\CC^2$ strictly
plurisubharmonic function on $M$. Assume that the set
$Diff_{(J,J_{st})}(D,\mathbb B_2)$ is empty. Then the set
$Diff_{(J,J)}(D,D)$ is compact for the compact-open topology.
\end{theorem}

\section{Preliminaries}
The following Lemma shows that every almost complex structure $J$ on an
almost complex manifold may be represented locally at $p \in M$
 as a small $\mathcal C^2$ deformation of the standard structure
(see~\cite{gs03}). 

\begin{lemma}
\label{suplem1}
For every $\lambda_0 > 0$ there exist a neighborhood $U_0$ of $p$ and a
coordinate diffeomorphism 
$z : U_0 \rightarrow \mathbb B_2$ such that
$z(p) = 0$, $dz(p) \circ J(p) \circ dz^{-1}(0) = J_{st}$  and the
direct image $\hat J = z_*(J)$ satisfies $\vert\vert \hat J - J_{st}
\vert\vert_{\CC^2(\bar{\mathbb B}_2)} \leq \lambda_0$.
\end{lemma}
\proof Shrinking $U$ if necessary there exists a diffeomorphism $z$ 
from $U$ onto $\mathbb B_2$ satisfying $z(p) = 0$ and $dz(p) \circ J(p)
\circ dz^{-1}(0) = J_{st}$. For $\lambda > 0$ consider the dilation
$d_{\lambda} : t \mapsto \lambda^{-1}t$ in $\C^2$ and the composition
$z_{\lambda} = d_{\lambda} \circ z$. Then $\lim_{\lambda \rightarrow
0} \vert\vert (z_{\lambda})_{*}(J) - J_{st} \vert\vert_{\CC^2(\bar
{\mathbb B}_2)} = 0$. Setting $U = z^{-1}_{\lambda}(\mathbb B_2)$ for
$\lambda > 0$ small enough, we obtain the desired statement. \qed

\vskip 0,1cm
In the sequel we will consider the diffeomorphism $z$ given by
Lemma~\ref{suplem1}, for sufficiently small~$\lambda_0$.

\subsection{Plurisubharmonic functions}
Let $(M,J)$ be an almost complex manifold. We denote by $TM$ the real 
tangent bundle of $M$ and by $T_\C M$ its complexification. If
$T^{(1,0)}M :=\{ X \in T_\C M  : JX=iX\} = \{\zeta -iJ \zeta, \zeta \in TM\}$,
$T^{(0,1)}M :=\{ X \in T_\C M : JX=-iX\} = \{\zeta +iJ \zeta, \zeta \in TM\}$
then $T_\C M = T^{(1,0)}M \oplus T^{(0,1)}M$. The set of complex
forms of type $(1,0)$ on $M$ is defined by
$T_{(1,0)}M=\{w \in T_\C^* M : w(X) = 0, \forall X \in T^{(0,1)}M\}$
and the set of complex forms of type $(0,1)$ on $M$ by
$T_{(0,1)}M=\{w \in T_\C^* M : w(X) = 0, \forall X \in T^{(1,0)}M\}$.
Then $T_\C^*M=T_{(1,0)}M \oplus T_{(0,1)}M$ and the operators $\partial_J$ and
$\bar{\partial}_J$ are defined on the space of smooth complex functions on
$M$ by $\partial_J u =
du_{(1,0)} \in T_{(1,0)}M$ and $\bar{\partial}_Ju = du_{(0,1)} \in T_{(0,1)}M$,
for every complex smooth function $u$ on $M$. 

\vskip 0,1cm
Let $r$ be a $\CC^2$ function on $M$. The Levi 
form of $r$ is defined on $TM$ by 
$$
\mathcal L^J(r)(X) :=-d(J^\star dr)(X,JX).
$$
We recall that an upper semicontinuous function $u$ on $(M,J)$ is called 
{\sl $J$-plurisubharmonic} on $M$ if the composition $u \circ f$ 
is subharmonic on $\Delta$ for every $f \in \mathcal O_J(\Delta,M)$.
Then we have the following characterization of
$J$-plurisubharmonic functions (see \cite{gs03}) :
\begin{proposition}\label{PROP}
Let $u$ be a $\CC^2$ real valued function on $M$. Then $u$ is
$J$-plurisubharmonic on $M$ if and only if $\mathcal L^J(u)(X) \geq 0$
for every $X \in TM$.
\end{proposition}

\vskip 0,1cm
Hence, following \cite{gs03}, we say that a $\CC^2$ real valued function 
$u$ on $M$ is {\sl strictly $J$-plurisubharmonic} on $M$ if 
$\mathcal L^J(u)$ is positive definite on $TM$.

We point out that the strict $J$-plurisubharmonicity is stable with respect to
small $\mathcal C^2$ deformations of the almost complex structure and of the
function.
\subsection{Kobayashi-Royden pseudometric}
We denote by $\Delta$ the unit disc in $\C$. Let $(M,J)$ be an almost
complex manifold. According to \cite{nw63} for
every $p \in M$ there is a neighborhood $\mathcal V$ of $0$ in $T_pM$ such
that for every $v \in \mathcal V$ there exists $f \in \mathcal
O_J(\Delta,M)$ satisfying $f(0) = p,$ $df(0) (\partial / \partial x) =
v$. We may therefore define the Kobayashi-Royden infinitesimal pseudometric
$K_{(M,J)}$ and the integrated pseudodistance $d^K_{(M,J)}$ (the upper
semicontinuity of $K_{(M,J)}$ on the tangent bundle $TM$ of $M$ is proved in
\cite{kr99})~:
\begin{definition}\label{dd}
$(i)$ For $p \in M$ and $v \in T_pM$, $K_{(M,J)}(p,v)$ is the infimum of the
set of positive $\alpha$ such that there exists a $J$-holomorphic disc
$f :\Delta \rightarrow M$ satisfying $f(0) = p$ and $df(0)(\partial
/\partial x) = v/\alpha$.

$(ii)$ Let $p,q \in M$. Denote by $\Gamma_{p,q}$ the set of all $\mathcal
C^1$-paths $\gamma :[0,1] \rightarrow M$ satisfying $\gamma(0) =p,\,
\gamma(1) = q$. Then $
d_{(M,J)}^K(p,q)=\inf_{\gamma \in \Gamma_{p,q}}\int_0^1K_{(D,J)}
(\gamma(t),\gamma'(t))dt$.
\end{definition}

As in the complex case, the integrated pseudodistance is decreasing under
the action of $(J',J)$-holomorphic maps.
Hence if $f \in Diff_{(J',J)}(M',M)$ then the inverse map
$f^{-1}$ is in 
$Diff_{(J,J')}(M,M')$ and we have for every $p',q' \in M'$~: 
\begin{equation}\label{equation4}
d_{(M,J)}^K(f(p'),f(q')) = d_{(M',J')}^K(p',q').
\end{equation}
\vskip 0,1cm
The following results are proved in \cite{gs03} :

\vskip 0,1cm
\noindent{\bf Proposition A.} {\it (Localization principle)
Let $D$ be a domain in an almost complex
manifold $(M,J)$, let $p \in \bar{D}$ and let $U$ be a neighborhood of
$p$ in $M$ (not necessarily contained in $D$). 
Let $u$ be a $\mathcal C^2$ function on $\bar{D}$, negative and
$J$-plurisubharmonic on $D$. We assume that $-L \leq u < 0$ on $D \cap
U$ and that $u - c |z|^2$ is strictly $J$-plurisubharmonic on $D \cap U$,
where $c$ and $L$ are positive constants and $z$ is the diffeomorphism
given by Lemma~\ref{suplem1}. Then there exist a positive
constant $ s$ and a neighborhood $V \subset \subset U$ of $p$,
depending on $c$ and $L$ only, such that for $q \in D \cap V$ and $v
\in T_qM$ we have $K_{(D,J)}(q,v) \geq s K_{(D \cap U,J)}(q,v)$.
}

\vskip 0,1cm
The next Proposition gives a lower bound on the Kobayashi-Royden infinitesimal
pseudometric (more precise lower estimates are given in Theorem~1 of
\cite{gs03}).

\vskip 0,1cm
\noindent{\bf Proposition B.} {\it Let $D$ be a relatively compact
domain in an almost complex manifold $(M,J)$.
Let $u$ be a $\mathcal C^2$ function on $\bar{D}$, satisfying
$-L \leq u < 0$ on $D$ and $u - c|z|^2$ is strictly
$J$-plurisubharmonic on $D$, where
$c$ and $L$ are positive constants and $z$ is given by Lemma~\ref{suplem1}.
Then there exist positive constants $C$ and
$\lambda_0$, depending on $c$ and $L$ only,
such that for every almost complex structure $J'$ defined
in a neighborhood of $\bar{D}$ and such that
$\|J'-J\|_{\CC^2(\bar{D})} \leq \lambda_0$ we have~:
$K_{(D,J')}(p,v) \geq C \|v\|$
for every $p \in D$ and every $v \in T_pM$.}

\vskip 0,1cm
Finally we have the boundary behaviour of the Kobayashi pseudodistance
in a strictly $J$-pseudoconvex domain.

\vskip 0,1cm
\noindent{\bf Proposition C.} {\it Let $D$, $D^\nu$ be domains in an
almost complex manifold $(M,J)$ and let $p \in \partial D$.
Assume that there is a neighborhood $U$ of $p$ such that 
$D \cap U =\{x \in U :r(x)<0 \}$, where $r$ is a $\CC^2$ 
strictly $J$-plurisubharmonic function in a neighborhood of $\bar{U}$. 
If $D^\nu \cap U = \{x \in U : r^\nu(x) < 0 \}$, where
$r^\nu$ is a sequence of $\mathcal C^2$ functions in a neighborhood
of $\bar{U}$, converging to $r$
in the $\mathcal C^2(\bar{U})$ convergence, then for every $q \in D \cap U$
there is $\nu_0$ such that for $\nu \geq \nu_0$ we have~:
$\lim_{\stackrel{x \in D^\nu}{x \rightarrow p}}d_{(D^\nu,J)}^K(x,q) =
+\infty$.}


\section{Proof of theorem~\ref{thm2.1}}
In this Section we assume that the assumptions of Theorem~\ref{thm2.1} are
satisfied. 

\subsection{Attraction property}
The following Lemma is a direct application of Proposition~C.
\begin{lemma}\label{lem3.3.1}
For every $K \subset \subset D$ we have : 
$\lim_{\nu \rightarrow \infty}f^\nu(K) =p$.
\end{lemma}
\noindent{\it Proof of Lemma \ref{lem3.3.1}}. Let $K \subset \subset D$ 
be such that $x^0 \in
K$. Since the function $x \mapsto d_D^K(x^0,x)$ is bounded from above by a
constant $C$ on $K$, it follows from the decreasing property of the Kobayashi 
pseudodistance that

\begin{equation}\label{eq2}
d_{(G^\nu,J_\nu)}^K(f^\nu(x^0),f^\nu(x)) \leq C
\end{equation}
for every $\nu$ and every 
$x \in K$. Moreover, $r$ is strictly $J$-plurisubharmonic 
in a neighborhood of 
$\bar{U}$ and the sequence $(r_\nu)_\nu$ converges to $r$ in the 
$\CC^2(\bar{U})$ convergence. It follows from Proposition~C that for
 every $V \subset \subset U$, containing $p$, we have :
\begin{equation}\label{eq3}
\lim_{\nu \rightarrow \infty}d_{(G^\nu,J_\nu)}^K 
(f^\nu(x^0),G^\nu \cap \partial V) = +\infty.
\end{equation} 
It follows from conditions (\ref{eq2}) and (\ref{eq3}) that 
$f^\nu(K) \subset V$ for every sufficiently large $\nu$. 
This gives the statement. \qed 

\vskip 0,1cm
According to \cite{sik94} Corollary~3.1.2,
there exist a neighborhood $U$ of $p$ in $M$ and complex coordinates
$z=(z_1,z_2)  : U \rightarrow \mathbb B_2 \subset \mathbb C^2$, $z(p) =
0$ such that $z_*(J)(0) = J_{st}$ and moreover, a map $f : \Delta
\rightarrow \mathbb B_2$ is $J' := z_*(J)$-holomorphic if it satisfies the
equations 
\begin{eqnarray}
\label{Jhol}
{\partial f_j}/ {\partial \bar \zeta} =
A_j(f_1,f_2)\overline{({\partial f_j}/{\partial \zeta}) },\ \ {\rm for}\ j=1,2,
\end{eqnarray}  
where $A_j(z) =  O(\vert
z \vert)$, $j=1,2$.

In order to obtain such coordinates, one can consider two
transversal foliations of $\mathbb B$ by $J'$-holomorphic curves (see
\cite{nw63}) and then take these curves into the lines $z_j = const$ by a local
diffeomorphism. The direct image of the almost complex structure $J$ under
such a diffeomorphism has a diagonal matrix $
J'(z_1,z_2) = (a_{jk}(z))_{jk}$ with $a_{12}=a_{21}=0$ and
$a_{jj}=i+\alpha_{jj}$ where $\alpha_{jj}(z)=\mathcal O(|z|)$ for $j=1,2$.

In what follows we omit the prime and denote this structure again by
$J$. We also still denote by $J_\nu$ the structure $z_*(J_\nu)$ defined
in $\mathbb B_2$ and we set $p^\nu :=z(f^\nu(x^0))$.
Hence we may assume that $U=\mathbb B_2$ and that $G^\nu \cap U=G^\nu \cap
\mathbb B_2$ are domains in $\mathbb C^2$, contained in $\mathbb B_2$.
Let $G :=\{x \in \mathbb B_2 : r(x) < 0\}$.
We may assume that the complex tangent space $T_0(\partial G)
\cap J(0) T_0(\partial G) = T_0(\partial G) \cap i T_0(\partial G)$ is
given by $\{ z_2 = 0 \}$.
In particular, we have the following expansion for the defining
function $r$ of $G$ on $\mathbb B_2$ :
$r(z,\bar{z}) = 2 Re(z_2) + 2Re K(z) + H(z,\bar{z}) + \mathcal O(\vert z
\vert^3)$, where
$K(z)  = \sum \lambda_{k,l} z_{k}{z}_{l}$, $\lambda_{k,l} =
\lambda_{l,k}$ and 
$H(z,\bar{z}) = \sum \lambda_{k,\bar{l}} z_{k}\bar z_{l}$,
$\lambda_{k,\bar{l}} = \bar \lambda_{l,\bar{k}}$.

\begin{lemma}\label{PP}
The domain $G$  is strictly $J_{st}$-pseudoconvex near the origin.
\end{lemma}

\noindent{\sl Proof of Lemma~\ref{PP}}. Consider a complex
 vector $v=(v_1,0)$ tangent to $\partial D$ at the origin.
Let $f :\Delta \rightarrow \mathbb C^2$ be a  $J$-holomorphic disc 
centered at the origin and tangent to $v$ : $f(\zeta) = v\zeta +
 \mathcal O(\vert \zeta \vert^2)$.
Since $A_2 = \mathcal O(\vert z \vert)$, it follows from the $J$-holomorphy
equation (\ref{Jhol}) that 
$(f_2)_{\zeta\bar\zeta}(0) = 0$. This implies that 
$(r \circ f)_{\zeta\bar\zeta}(0) = H(v).$ Thus, the Levi form with
respect to $J$ coincides with the Levi form with respect to $J_{st}$
on the complex tangent space of $\partial G$ at the origin.\qed

\vskip 0,2cm
For sufficiently large $\nu$ let $q^\nu$ be the unique point on
$\partial G^\nu \cap \mathbb B_2$ such that 
$|q^\nu - p^\nu| = dist(p^\nu,\partial G^\nu \cap \mathbb B_2)$.
Since the sequence $(J_\nu)_\nu$ converges to $J$ in the 
$\CC^2(\bar{\mathbb B}_2)$ convergence,
there exists a diffeomorphism $A^\nu$ defined on
$\mathbb B_2$ such that $A^\nu(q^\nu)=0$ and if
$A^\nu(z)=:z^\nu=(z_1^\nu,z_2^\nu)$
then the expansion of $r_\nu$ in the $z^\nu$-coordinates is given by~:
$r_\nu = 2 \lambda_n^\nu Re(z_2^\nu)
+ 2 Re(\sum_{k,l =1}^2\lambda_{k,l}^\nu z_k^\nu z_l^\nu) 
+ \sum_{k,l=1}^2\lambda_{k,\bar{l}}^\nu z_k^\nu \bar{z}_l^\nu + 
\mathcal O(|z^\nu|^3)$, $\lambda_{k,l}^\nu = \lambda_{l,k}^\nu,\ 
\lambda_{k,\bar{l}}^\nu = \bar \lambda_{l,\bar{k}}^\nu$,
where the condition $\mathcal O(|z^\nu|^3)$ is
uniform with respect to $\nu$. Finally the diffeomorphism $A^\nu$ can
be chosen such that $d(A^\nu) \circ J_\nu \circ d(A^\nu)^{-1}$
is represented by a diagonal matrix $(a_{jk}^\nu(z^\nu))_{jk}$ with
$a_{12}^\nu=a_{21}^\nu=0$ and $a_{jj}^\nu = i + \alpha_{jj}^\nu(z^\nu)$ where
$a_{jj}^\nu$ converges to $a_{jj}$, with its first derivatives, uniformly
on compact subsets of $\mathbb B_2$. We note that $A^\nu$ converges
to the identity map in any $\mathcal C^k(\bar{\mathbb B}_2)$ norm. 
Since $r_\nu$ converges to $r$ in the $\mathcal C^2(\bar{\mathbb B}_2)$ norm
by assumption, it follows that~: $\lim_{\nu \rightarrow \infty}\lambda_2^\nu
= 1$, $lim_{\nu \rightarrow \infty}\lambda_{k,l}^\nu = \lambda_{k,l}$
and 
$\lim_{\nu \rightarrow \infty}\lambda_{k,\bar{l}}^\nu = \lambda_{k,\bar{l}}$
for $k,l=1,2$.

\vskip 0,2cm 
For every $\nu \geq 1,$ let $\tau_\nu := dist(p^\nu,q^\nu)$.
We define the dilation $\Lambda^\nu : A^\nu(\mathbb B_2) \rightarrow
\mathbb C^2$ by 
$\Lambda^\nu(z^\nu) = (z_1^\nu / \sqrt{\tau_\nu}, z_2^\nu / \tau_\nu)$,
in the $z^\nu$-coordinates.
Let $\tilde{J}_\nu$ be the almost complex structure defined on
$\Lambda^\nu \circ A^\nu(\mathbb B_2)$
by $\tilde{J}_\nu :=\Lambda^\nu \circ dA^\nu \circ J_\nu 
\circ d(A^\nu)^{-1} \circ (\Lambda^\nu)^{-1}$.
Since $A^\nu$ is $(J_\nu, dA^\nu \circ J_\nu 
\circ d(A^\nu)^{-1})$-holomorphic, $\Lambda^\nu$ is 
$(J_\nu,\tilde{J}_\nu)$-holomorphic. Finally, we consider the domain
$$
\tilde{G}^\nu :=\Lambda^\nu \circ A^\nu(G^\nu \cap \mathbb B_2) =
\{z \in \Lambda^\nu \circ A^\nu (\mathbb B_2) : 
r_\nu((\Lambda^\nu \circ A^\nu)^{-1}(z),
\overline{(\Lambda^\nu \circ A^\nu)^{-1}(z)}) <0\}.
$$
The next Lemma gives the limit behaviour of the domains $\tilde{G}^\nu$ 
and of the almost complex structures $\tilde{J}_\nu$.
\begin{lemma}\label{lem3.3.2}
The following conditions are satisfied :

$(i)$ $\lim_{\nu \rightarrow \infty}\tilde{G}^\nu = 
\mathbb
G=\{z \in \C^2 : Re(z_2) + 2 Re(K(z_1,0)) + H((z_1,0),(\bar{z}_1,0)) < 0\}$,

$(ii)$ $\lim_{\nu \rightarrow \infty}
\tilde{J}_\nu=J_{st}$ uniformly on compact subsets of $\mathbb C^2$.
\end{lemma}

\noindent{\it Proof of Lemma \ref{lem3.3.2}}. The convergence in
statement~$(i)$ is the Hausdorff convergence on compact
subsets. Condition~$(i)$ is a direct consequence of the
convergence of $(r_\nu)_\nu$ to $r$ and is similar to the usual complex case
(see~\cite{pi89}).

\noindent{\it Proof of $(ii)$}.
For every $\nu$ the almost complex structure $\tilde{J}_\nu$
is represented by the diagonal matrix $(b_{jk}^\nu(z))_{jk}$ where
$b_{jj}^\nu(z)=a_{jj}^\nu((\Lambda^\nu)^{-1}(z))
= i + \alpha_{jj}^\nu(\sqrt{\tau_\nu} z_1,\tau_\nu z_2)
\rightarrow_{\nu \rightarrow \infty}i$, uniformly on compact subsets of
$\mathbb C^2$. \qed

\vskip 0,1cm
For every $\nu$, let $F^\nu := \Lambda^\nu \circ A^\nu \circ f^\nu$.
Then $F^\nu$ is a $(J_{st},\tilde{J}_\nu)$-holomorphic map
from $(f^\nu)^{-1}(G^\nu \cap \mathbb B_2)$ to $\tilde{G}^\nu$ satisfying
$F^\nu(x^0) = (0,-1)$ and we have~:
\begin{proposition}\label{re-scaling}
$(i)$ We may extract from $(F^\nu)_\nu$ a subsequence converging, 
uniformly on compact subsets of $D$, to a $J_{st}$-holomorphic map 
$F : D \longrightarrow \bar{\mathbb G}$,

$(ii)$ We may extract from $((F^\nu)^{-1})_\nu$ a subsequence 
converging, uniformly on compact subsets of $\mathbb G$, to a
$J_{st}$-holomorphic map $\tilde{F} : \mathbb G \longrightarrow \bar{D}$,

$(iii)$ $F$ is a biholomorphism from $D$ to $\mathbb G$ with 
$F^{-1} = \tilde{F}$.
\end{proposition}
Statement~$(iii)$ is the conclusion of Theorem~\ref{thm2.1}. Indeed it follows
from Lemma~\ref{PP} that $H((z_1,0),(\bar{z}_1,0)) = \alpha z_1 \bar{z}_1$ with
$\alpha > 0$. Then the map 
$(z_1,z_2) \mapsto (\sqrt{\alpha} z_1,z_2 + H(z_1,0))$
is a biholomorphism from $\mathbb G$ to the unbounded representation
of the unit ball 
$\mathbb H = \{(z_1,z_2) \in \mathbb C^2 : Re(z_2) + |z_1|^2 < 0\}$.

\vskip 0,1cm
\noindent Our proof of Proposition~\ref{re-scaling} is based on the method
developped by F.Berteloot and G.Coeur\'e~\cite{be-co91} and by
F.Berteloot~\cite{be95}.
We first prove the following Lemma~:
\begin{lemma}\label{LEM}
There exist $C_0 > 0$, $\delta_0 > 0$ and $r_0 > 0$ such that for every
$0 < \delta < \delta_0$, for every $\nu >> 1$ and for every
${J}_\nu$-holomorphic disc $g^\nu : \Delta \rightarrow {G}^\nu$
we have~:
$$
g^\nu(0) \in Q(0,\delta) \Rightarrow g^\nu(r_0 \Delta) \subset Q(0,C_0 \delta)
,
$$
where $Q(0,\delta):=\{z \in \mathbb C^2 : |z_1|< \sqrt{\delta},
|z_2| < \delta\}$.
\end{lemma}
\noindent{\it Proof of Lemma~\ref{LEM}}. Assume by contradiction that there
exist $C_\nu \rightarrow \infty$, 
$\zeta_\nu \rightarrow 0$ in $\Delta$ and ${J}_\nu$-holomorphic discs
$g^\nu : \Delta \rightarrow {G}^\nu$ such that
$g^\nu(0) \in Q(0,\delta_\nu)$ 
and
$g^\nu(\zeta_\nu) \not\in Q(0,C_\nu \delta_\nu)$. 
Let $d_\nu : z \mapsto (z_1/\sqrt{\delta_\nu}, z_2/\delta_\nu)$, let
$h^\nu:=d_\nu \circ g^\nu$ and let $J^\nu:=d_\nu(J_\nu)$. We set
$H^\nu:=\{z \in d_\nu(U) : r_\nu \circ (d_\nu)^{-1}(z,\bar{z}) < 0\}$.
It follows from Lemma~\ref{lem3.3.2} part $(ii)$ that
$\rho^\nu:=r_\nu \circ d_\nu$ converges to $\rho:=Re(z_2) + |z_1|^2$,
uniformly on compact subsets of $\mathbb C^2$ and $J^\nu$ converges
to $J_{st}$, uniformly on compact subsets of $\mathbb C^2$.
There exists $A>0$ such that the function $\rho + A \rho^2$ is strictly
$J_{st}$-plurisubharmonic on $Q(0,5)$. Hence for sufficiently large $\nu$
the function $\rho^\nu + A (\rho^\nu)^2$ is strictly $J^\nu$-plurisubharmonic
on $Q(0,4)$. Moreover we can extend this function as a
$J^\nu$-plurisubharmonic function on $H^\nu$. In particular it follows from
Proposition~A and Proposition~B
that there is a positive constant $C$ such that
$K_{(H^\nu,J^\nu)}(z,v) \geq C \|v\|$ for every $z \in H^\nu \cap Q(0,3),\
v \in \C^2$. Therefore, there exists a
constant $C' > 0$ such that $\parallel dh_\nu(\zeta) \parallel \leq C'$ for
any $\zeta \in (1/2)\Delta$ satisfying $h_\nu(\zeta) \in H^\nu \cap Q(0,3)$. On
the other hand, the sequence $\vert h_\nu(\zeta_\nu) \vert$ tends to $+
\infty$. Denote by $[0,\zeta_\nu]$ the segment 
(in $\mathbb C$) joining the origin and $\zeta_\nu$ and let 
$\zeta_\nu' \in [0,\zeta_\nu]$ be the point closest to the origin such
that  $h_\nu([0,\zeta_\nu']) \subset H^\nu \cap
\overline Q(0,2)$ and $h_\nu(\zeta_\nu') \in \partial Q(0,2)$. Since $h_\nu(0)
\in Q(0,1)$, we have  $\vert h_\nu(0) - h_\nu(\zeta_\nu') \vert \geq C''$
for some constant $C'' > 0$. Let $\zeta_\nu' = r_\nu e^{i\theta_\nu}$,
$r_\nu \in ]0,1[$. Then 
$$
\vert h_\nu(0) - h_\nu(\zeta_\nu') \vert \leq \int_{0}^{r_\nu}
\parallel dh_\nu(te^{i\theta_\nu}) \parallel dt \leq C'r_\nu\longrightarrow 0.
$$
This contradiction proves Lemma~\ref{LEM}.

\vskip 0,1cm
We note that Lemma~\ref{LEM} is also satisfied replacing $J_\nu$-holomorphic
discs by $(J_{st},J_\nu)$-holomorphic maps. 
As a corollary we have the following
\begin{lemma}\label{conv}
For any compact subset $K \subset D$ the sequence of norms
$(\|F^\nu\|_{\mathcal C^0(K)})_\nu$ is bounded.
\end{lemma}
For the proof we can consider a covering of $K$
by sufficiently small balls, similarly to \cite{be-co91},
p.84. Indeed, consider a covering of $K$ by the balls $p^j + r_0 \mathbb B$,
$j=0,...,N$ where $r_0$ is given by Lemma~\ref{LEM} and
$p^{j+ 1} \in p^j + r_0 \mathbb B$ for
any $j$. Since $\lim_{\nu \rightarrow \infty} (A^\nu \circ f^\nu)(x^0) = 0$,
 we obtain, for $\nu$ large enough,
that $A^\nu \circ f^\nu(p^0 + r\mathbb B) \subset Q(0,2C\tau_\nu)$, then
$A^\nu \circ f^\nu(p^1 + r\mathbb B) \subset Q(0,4C^2\tau_\nu)$.
Continuing this process we obtain that
$A^\nu \circ f^\nu(p^N + r\mathbb B) \subset Q(0,2^NC^N\tau_\nu)$.
Since $\tau_\nu \rightarrow 0$ and
$F^\nu=\Lambda^\nu \circ A^\nu \circ f^\nu$ we obtain Lemma~\ref{conv}.

\vskip 0,1cm
\noindent We prove now Proposition~\ref{re-scaling}. {\it Part $(i)$}. 
Lemma~\ref{conv} implies that the sequence
$(F^\nu)_\nu$ is bounded (in the $\mathcal C^0$ norm) on any
compact subset $K$ of $D$. Covering $K$ by small bidiscs,
consider two transversal foliations by holomorphic curves on every
bidisc. Since the restriction of $F^\nu$ on every such curve is
uniformly bounded in the $C^0$-norm, 
it follows by the well-known elliptic
estimates that it is bounded in $C^l$ norm for every $l$ (see
\cite{sik94}). Since the bounds are uniform with respect to curves, this
implies that the sequence $(F^\nu)_\nu$ is bounded in every
$C^l$-norm. So the family $(F^\nu)_\nu$ is relatively compact by Ascoli
theorem. Let $F$ be a cluster point of $(F^\nu)_\nu$. We still denote by
$(F^\nu)_\nu$ an appropriate subsequence converging (uniformly on compact
subsets of $D$) to $F$.
Passing to the limit in the holomorphy
condition $\tilde{J}_\nu \circ dF^\nu = dF^\nu \circ \tilde{J}_\nu$,
we obtain that $F$ is holomorphic with respect to
 $J_{st}$.

\vskip 0,1cm
\noindent {\it Part $(ii)$}. 
Let $0 < r < 1$, $B(0,r) = \{z \in \C^2 : \|z\|<r\}$ and let $\Phi$ be a
biholomorphism from $\mathbb G$ to $\mathbb B_2$ satisfying
$\Phi((0,-1))=(0,0)$ (this is a biholomorphism
for the standard structure $J_{st}$ on $\mathbb G$ and $\mathbb B_2$). 
According to Lemma~\ref{lem3.3.2} $(ii)$
we have, for sufficiently large $\nu$, the inclusion 
$\Phi^{-1}(B(0,r)) \subset \tilde{G}^\nu$.
Since the sequence 
$(\tilde{J}_\nu)_\nu$ converges uniformly on compact subsets of $\mathbb G$ 
to $J_{st}$ (in $\mathcal C^2$ norm) by
Lemma~\ref{lem3.3.2} $(iii)$, the sequence 
$(J'_\nu :=d\Phi \circ \tilde{J}_\nu \circ d\Phi^{-1})_\nu$ converges to
$J_{st}$ uniformly on compact subsets of $\mathbb B_2$
(in $\mathcal C^2$ norm). 
Fix $0<r<1$. Then there exists a positive constant $C_r$ and $r < r' < 1$
such that for every $z \in B(0,r)$ and every unitary vector $v$ in
$\mathbb C^2$ there is a $J_{st}$-holomorphic disc $\varphi_{z,v} :
\Delta \rightarrow \mathbb B_2$ satisfying $\varphi_{z,v}(0) = z,\
\varphi'(0) \in \mathbb Cv$, $|\varphi'(0)| \geq C_r$ and
$\varphi_{z,v}(\Delta(0,1/2)) \subset B(0,r')$.
According to Section~5.4a of 
\cite{nw63} (study of the stability of $J$-holomorphic discs under small
smooth deformations of a given almost complex structure $J$)
there exists a positive contant $C'_r$ and $r < r'' < 1$
such that for sufficiently large $\nu$, for every 
$z \in B(0,r)$ and for every unitary vector
$v\in \C^2$, there is a $J'_\nu$-holomorphic disc 
$f_{J'_\nu,z,v} : \Delta \rightarrow B_2$, centered at $z$, satisfying 
\begin{equation}\label{EQUA}
d(f_{J'_\nu,z,v})(0)(\partial / \partial x) = \alpha v\ \ 
{\rm with} \ \alpha \geq C'_r
\end{equation}
and $f_{J'_\nu,z,v})(\Delta(0,1/3)) \subset B(0,r'')$.
Since for sufficiently large $\nu$ the $J_{st}$-holomorphic disc
$\tilde{f}_{J'_\nu,z} :=(\Phi \circ F^\nu)^{-1} \circ f_{J'_\nu,z}$ 
satisfies the inclusion 
$\tilde{f}_{J'_\nu,z}(\Delta) \subset D \subset \subset \C^2$,
there exists a positive constant $C'$ such that
\begin{equation}\label{Equa2}
|(\tilde{f}_{J'_\nu,z})'(0)| \leq C'
\end{equation}
for sufficiently large $\nu$.

It follows from conditions~(\ref{EQUA}) and (\ref{Equa2}) that the first 
derivative of $(\Phi \circ F^\nu)^{-1}$ is bounded from below by a positive 
constant on $B(0,r)$, uniformly with respect to $\nu > > 1$.
By Ascoli Theorem we may extract from $((F^\nu)^{-1})_\nu$ a 
subsequence that converges uniformly on $B(0,r)$ to a holomorphic map 
$\tilde F_r : B(0,r) \rightarrow \bar{D}$. 
In particular $\tilde F_{r'} = \tilde F_r$ on $B(0,r)$ for $r' > r$.
This proves condition $(ii)$.

\noindent{\it Part $(iii)$}. 
We know that $F(x^0) = 0 \in \mathbb G$. Assume now that there is $x
\in D$ such that $F(x) \in \partial \mathbb G$ and let $\gamma$ be a
$\mathcal C^1$-path in $D$ such that $\gamma(0) = x^0$, $\gamma(1) =
x$. We consider $t_0 \leq 1$ such that $F(\gamma(0,t_0[)) \subset
\mathbb G$ and $F(\gamma(t_0)) \in \partial \mathbb G$. Since $\mathbb
G$ is complete hyperbolic ($\mathbb G$ is biholomorphic to the unit
ball in $\C^2$), we obtain that $\lim_{t \rightarrow t_0}d_\mathbb
H^K(0,F(\gamma(t)))= \infty$. However, for every $t <1$~:
$$
d_\mathbb G^K(0,F(\gamma(t))) \leq 
\sup_{s \in [0,1]}d_D^K(x^0,\gamma(s))<\infty$$
by the compactness of $\gamma([0,1])$ in $D$. This is a contradiction,
 implying that $F(D) \subset \mathbb G$. We prove now that $\tilde{F}(\mathbb
 G) \subset D$. Since $D$ is bounded in $\C^2$, $D$ is hyperbolic for
 the standard structure $J_{st}$. Let $r>0$ and consider the Kobayashi
 ball $B_{(D,J_{st})}^K(x^0,r)$. Since $f^\nu$ is a biholomorphism from
 $(D,J_{st})$ to $(G^\nu,J_\nu)$, we have
 $f^\nu(B_{(D,J_{st})}^K(x^0,r))=B_{(G^\nu,J_\nu)}(f^\nu(x^0),r)$
(see equality~(\ref{equation4})). Since the sequence of points
 $(f^\nu(x^0))_\nu$ converges to $p$ we obtain from Proposition~C that
 $B_{(G^\nu,J_\nu)}(f^\nu(x^0),r) \subset \subset G^\nu \cap U$ for
 sufficiently large $\nu$. In particular, since $(f^\nu)^{-1}$ is
 continuous on $G^\nu$, we have : $B_{(D,J_{st})}^K(x^0,r) \subset
 \subset D$. This implies that $D$ is complete hyperbolic. Assume that
 $q \in \mathbb G$ is such that $\tilde{F}(q) \in \partial D$. Then there
 exists $r>0$ such that $q \in B_{(\mathbb G,J_{st})}^K( 0,r)$. Let
 $\gamma :[0,1] \rightarrow \mathbb G$ be a $\mathcal C^1$-path such
 that $\gamma(0) =0$, $\gamma(1) = q$ and $\gamma([0,1]) \subset
 B_{(\mathbb G,J_{st})}^K(0,2r)$. We define $t_0 \leq 1$ such that
 $\tilde{F}(\gamma([0,t_0[)) \subset D$ and 
$\tilde{F}(\gamma(t_0)) \in \partial
 D$. By the complete hyperbolicity of $D$, we have : $\lim_{t
 \rightarrow t_0}d_{(D,J_{st})}^K(x^0,\tilde{F}(\gamma(t))) = \infty$, which
 contradicts the condition $d_{(D,J_{st})}^K(x^0,\tilde{F}(\gamma(t))) \leq
 d_{(\mathbb G,J_{st})}(0,\gamma_t) \leq 2r$ for every $t <
 t_0$. Consequently : $\tilde{F}(\mathbb G) \subset D$. Now, since for every
 $\nu$, $(F^\nu)^{-1} \circ F^\nu= id_D$ and $F^\nu\circ (F^\nu)^{-1}
 = id_\mathbb G$ we obtain, by taking the limit when $\nu \rightarrow
 \infty$, that $\tilde{F} \circ F=id_D$ and 
$F \circ \tilde{F} = id_\mathbb G$, 
meaning that $F$ is a biholomorphism from $D$ to $\mathbb G$. 
This proves Theorem~\ref{thm2.1}. \qed

\vskip 0,1cm
Corollary~\ref{cor2.2} is a direct consequence of Theorem~\ref{thm2.1}.
 
\section{Proof of Theorem~\ref{thm2.3}}
We explain now how to adapt the proof of Theorem~\ref{thm2.1} to obtain
Theorem~\ref{thm2.3}. Hence we assume in this Section that the assumptions
of Theorem~\ref{thm2.3} are satisfied. 

\vskip 0,1cm
We have the following
\begin{lemma}\label{lem4.1}
For every $0 \in K \subset \subset \Delta$ we have~:
$\lim_{\nu \rightarrow \infty}(f^\nu \circ f)(K) =p$,
uniformly with respect to $f \in \mathcal O_{(J_{st},J)}(\Delta,D)$ such that
$f(0) = x^0$.
\end{lemma}

\noindent{\it Proof of Lemma~\ref{lem4.1}}. The proof follows line by line 
the proof of Lemma~\ref{lem3.3.1}, replacing $f^\nu$ by $f^\nu \circ f$ where
$f \in \mathcal O_{(J_{st},J)}(\Delta,D)$, $f(0) = x^0$. Indeed since the 
function $x \mapsto d_{(\Delta,J_{st})}^K(0,x)$ is bounded from above on $K$
by a positive constant $C$, we have~:
$$
d_{(D,J)}^K(f^\nu \circ f(0),f^\nu \circ f(x)) \leq C
$$
for every $\nu \geq 1$ and for every $x \in K$.
Fix a neighborhood $V \subset \subset M$ of $p$. Since 
$$
\lim_{\nu \rightarrow \infty}d_{(D,J)}^K(f^\nu(x^0),D \cap \partial V) 
= \infty
$$
we have $f^\nu \circ f(K) \subset V$
for sufficiently large $\nu$, uniformly with respect to 
$f \in \mathcal O_{(J_{st},J)}(\delta,D)$ such that $f(0) = x^0$. \qed

\vskip 0,1cm
\noindent{\it Proof of Theorem~\ref{thm2.3}}. Assume by contradiction that 
there is a sequence $(f^\nu)_\nu$ in $Diff_{(J,J)}(D,D)$ and  points
$x^0 \in D,\ p \in \partial D$ such that
$\lim_{\nu \rightarrow \infty}f^\nu(x^0)=p$.
Lemma \ref{lem4.1} implies that~: $\lim_{\nu \rightarrow \infty}f^\nu(K) = p$,
for every compact $K$ in $D$.

Let $U$, $z$ and $A^\nu$
defined as in page 5 and let $\Lambda^\nu$ is the dilation defined page 6.
We consider $G^\nu :=D$, $J_\nu :=J$ and $\tilde{J}_\nu :=\Lambda^\nu
\circ d(A^\nu) \circ J \circ d(A^\nu)^{-1} \circ (\Lambda^\nu)^{-1}$ 
for every $\nu$.
We may assume as in the proof of Theorem~\ref{thm2.1} that $D \cap U$ is in
$\mathbb C^2$ and we still denote by $J$
the associated structure in $\mathbb C^2$. If we set
$\tilde{G}^\nu :=\Lambda^\nu \circ A^\nu(D \cap U)$
then by Lemma~\ref{lem3.3.2}
we have $\lim_{\nu \rightarrow \infty}\tilde{G}^\nu = \mathbb G$
and
$\lim_{\nu \rightarrow \infty}\tilde{J}_\nu = J_{st}$,
uniformly on compact subsets of $\C^2$. If
 $F^\nu :=\Lambda^\nu \circ A^\nu \circ f^\nu$ 
then 
$F^\nu \in Diff_{(J,\tilde{J}_\nu)}((f^\nu)^{-1}(D \cap U), \tilde{G}^\nu)$
and according to Proposition~\ref{re-scaling} the sequence $(F^\nu)_\nu$
converges after extraction to a $(J,J_{st})$-holomorphic map $F$ from
$D$ to $\mathbb G$.
\vskip 0,1cm
We use the following quantitative version of Proposition~2.3.6 of 
\cite{sik94} :
\begin{proposition}\label{Prop}
Let $D'$ be a domain in $\C^n$. There is a positive 
constant $\delta_0$ such that for every almost complex structure
$J'$ in a neighborhood of $\overline{D'}$ satisfying
$\|J'-J_{st}\|_{\mathcal C^2(\overline{D'})} \leq \delta_0$
we have 
\begin{equation}\label{ineq1}
\|f_{J'}\|_{\mathcal C^1(K)} \leq C_{K,\delta_0}\|f_{J'}\|_{\mathcal C^0(K)},
\end{equation}
for every $f_{J'} \in \mathcal O_{(J_{st},J')}(\Delta,D')$ and for
every $K \subset \subset \Delta$, where $C_{K,\delta_0}$ is a positive
constant depending only on $K$ and $\delta_0$.
\end{proposition}
Let $U' \subset \subset U$ be a neighborhood of $p$ such that 
\begin{equation}\label{EQUATION}
\|J - J_{st}\|_{\mathcal C^2(\overline{U'})} \leq
\delta_0.
\end{equation}
Fix $0 < r < 1$, sufficiently close to 1.
Since $\lim_{\nu \rightarrow \infty}
\|\tilde{J}_\nu - J_{st}\|_{\mathcal C^2(\bar{B}(0,r))} = 0$, it follows from
\cite{nw63} Section~4.5a that there is a covering $\mathcal R_\nu$ of 
$\bar{B}(0,r)$ by $\tilde{J}_\nu$-holomorphic discs centered at the origin
for sufficiently large $\nu$, these discs being small deformations of the
straight holomorphic discs in the ball.
More precisely there exists $0 \in K_r \subset \subset \Delta$ and there
exists a positive constant $c_r$ such that for sufficiently large $\nu$
we have
$\bar{B}(0,r) \subset \displaystyle \cup_{f \in \mathcal R_\nu}f(K_r)$
and
\begin{equation}\label{Equation4.2}
\inf_{f \in \mathcal R_\nu, \zeta \in K_r}\|df(\zeta)(\partial / \partial x)\|
\geq c_r.
\end{equation}
Moreover we may assume that $\Phi^{-1}(\bar{B}(0,r)) \subset G^\nu$.
For $f \in \mathcal R_\nu$ consider the $(J_{st},J)$-holomorphic map
$g^\nu :=(F^\nu)^{-1} \circ \Phi^{-1} \circ f$
from $\Delta$ to $D$. Since $g^\nu(0) = x^0$ it follows from
Lemma~\ref{lem4.1}that there exists $\nu_0 \geq 1$ such that 
$f_{\nu_0} \circ g^\nu(K_r) \subset U'$,
uniformly with respect to $\nu >>1$.
It follows now from condition~(\ref{EQUATION}) and from Proposition~\ref{Prop}
($D \cap U' \subset \mathbb C^2$) that there exists a positive constant $C_r$
such that $\|f_{\nu_0} \circ g^\nu\|_{\mathcal C^1(K_r)} \leq C_r$
for sufficiently large $\nu$, 
or equivalently that
$\|g^\nu\|_{\mathcal C^1(K_r)} \leq C'_r$
with $C'_r >0$.

It follows from inequality~(\ref{Equation4.2}) that
$\|(F^\nu)^{-1} \circ \Phi^{-1}\|_{\mathcal C^1(\bar{B}(0,r))} \leq C''_r$
for a positive constant $C''_r$
and finally that
$\|(F^\nu)^{-1}\|_{\mathcal C^1(\Phi^{-1}(\bar{B}(0,r)))} \leq \tilde{C}_r$
for every sufficiently large $\nu$, where $\tilde{C}_r$ is a 
positive constant. 

Since $D$ is relatively compact in $M$, Ascoli Theorem implies that some
subsequence of $((F^\nu)^{-1})_\nu$ converges to a $(J_{st},J)$-holomorphic
map $\tilde{F}$ from $\mathbb G$ to $\bar{D}$. The end of the 
proof of Theorem~\ref{thm2.3} follows line by line the proof of 
Proposition~\ref{re-scaling} $(iii)$. Since $\tilde{F} \circ F = id_D$ and
$F \circ \tilde{F} = id_{\mathbb G}$ we obtain that 
$F \in Diff_{(J,J_{st})}(D,\mathbb G)$. 
This gives the contradiction. \qed

\vskip 1cm
{\small
\begin{tabular}{lllll}
{\sc Herv\'e Gaussier} & & & & {\sc Alexandre Sukhov}\\
C.M.I. & & & & U.S.T.L. \\
39, rue Joliot-Curie,& & & & Cit\'e Scientifique \\
13453 Marseille Cedex 13 & & & & 59655 Villeneuve d'Ascq Cedex\\
 & & & & \\
gaussier@cmi.univ-mrs.fr & & & & sukhov@agat.univ-lille1.fr
\end{tabular}
}
\end{document}